\documentclass[11pt]{article}
\usepackage{amssymb}
\newtheorem{lemma}{Lemma}[section]

\newtheorem{theorem}[lemma]{Theorem}
\newtheorem{proposition}[lemma]{Proposition}
\newtheorem{corollary}[lemma]{Corollary}
\newtheorem{example}[lemma]{Example}

\newtheorem{remark}[lemma]{Remark}

\def\endproof{\hfill$\Box$}

\def\endproof{\hfill$\Box$}

\title{Expansions and restrictions of structures, their hierarchies\footnote{The work
was carried out in the framework of Russian Scientific Foundation,
Project No.~24-21-00096.}}
\author{S.V. Sudoplatov}
\date{}

\begin{document}

\maketitle
\begin{abstract}
We introduce and study some general principles and hierarchical
properties of expansions and restrictions of structures and their
theories. The general approach is applied to describe these
properties for classes of $\omega$-categorical theories and
structures, Ehrenfeucht theories and their models, strongly
minimal, $\omega_1$-theories, and stable ones.
\end{abstract}


\bigskip
\section{Introduction}

The study and description of the possibilities of distributions of
expansions and restrictions of structures and their theories is of
interest both in general and for various natural algebraic,
geometric, ordered theories and models. The prerequisites for the
description are the original operations of Morleyzation, or
Atomization \cite{Hodges, Skolem}, and Skolemization \cite{ErPa,
Hodges, Skolem}, allowing to preserve or non-radically include
formulaically defined sets of a given structure, and upon receipt
of this elimination of quantifiers, according to which
formulaically defined levels arise in the form of Boolean
combinations of definable sets specified by quantifier-free
formulae. The operations of Shelahizations \cite{BeHr, Simon}, or
Namizations, produce both extensions and expansions of a structure
giving names or labels for definable sets. In addition, main types
of combinations of theories based on $E$-combinations and
$P$-combinations of their models are used for various expansions
preserving a series of properties \cite{cs}. These data allow to
introduce topological approximations and develop methods of
approximation \cite{at}, closures \cite{cl}, rank characteristics
of a family of theories \cite{RSrank} introduced and studied
within the framework of approximating formulae \cite{af} and, in
particular, pseudo-finite formulas \cite{pf}. In addition, the
effects associated with the preservation and violation of
properties during expansions and restrictions, with significant
formulaicity or non-formulaicity of these properties is
investigated. The formulaic properties include known theories with
the stability property \cite{Sh}, according to which a theory is
stable if each its formula is stable, i.e. it does not have the
order property, and a non-formulaic property, for example, is the
property of the theory of normality \cite{Pi83}, which, as the
author's examples show \cite{SudDis90}, may not be preserved when
taking restrictions. Other formulaic properties include the
properties of total transcendence and superstability \cite{Sh}
based on the monotonicity of ranks when passing to restricted
theories. A semantic characteristic is the value of the number of
non-isomorphic countable models of these theories, which can both
decrease and increase when the theory is restricted. Using these
tools, it is possible to describe the conditions for the
preservation and violation of significant properties when
constructing theories and their models, both in general and for
series of known classes of algebraic, geometric, ordered theories
and their models.

The paper is organized as follows. In Section 2, we introduce the
notion of regular structure, Boolean algebras for regular
expansions and restrictions based on a given universe both for
structures and their theories, and characterize these Boolean
algebras up to an isomorphism in terms of cardinalities of
universes. Kinds of properties on these Boolean algebras with
respect to upper and lower cones, lattices, and permutations are
studied in Sections 3, 4, 5, respectively. We apply these
divisions of properties to the classes of $\omega$-categorical
theories and structures, Ehrenfeucht theories and structures
(Section 6), strongly minimal ones (Section 7),
$\omega_1$-theories and structures (Section 8), stable ones
(Section 9). In particular, it is shown that any fusions of
strongly minimal structures are strongly minimal, too (Theorem
\ref{th_fus_sm}), whereas the properties of $\omega$-categoricity,
Ehrenfeuchtness, $\omega_1$-categoricity, and stability can fail
under fusions.

\section{Regular structures, their expansions, restrictions \\ and
Boolean algebras}

For convenience we consider {\em regular} structures, i.e.
relational structures without repetitions of interpretations of
signature symbols. There is no loss of generality with this
assumption, since operations can be replaced by their graphs, and
multiple names for the same relations are reduced to one of them.
The procedure transforming an arbitrary structure $\mathcal{M}$ to
a regular one $\mathcal{N}$ is called the {\em regularization},
and the structure $\mathcal{N}$ is called {\em regularized}, with
respect to $\mathcal{M}$. And the converse procedure transforming
$\mathcal{N}$ to the initial $\mathcal{M}$ is called {\em
deregularization}, and $\mathcal{M}$ is called {\em
deregularized}, with respect to $\mathcal{N}$.

\begin{remark}\label{reg_dereg}\rm
Clearly, regularizations and deregularizations are multi-valued
with respect to signatures, in general, but they always preserve
families of definable sets on the given universe.
\end{remark}

Let $\mathcal{M}$ be a regular structure. Let
$\overline{\mathcal{M}}$ be a maximal regular expansion of
$\mathcal{M}$ preserving the universe $M$. We have
\begin{equation}\label{eq_card_Sigma}
\left|\Sigma(\overline{\mathcal{M}})\right|={\rm
max}\{2^{|M|},\omega\},
\end{equation}
where $\Sigma(\overline{\mathcal{M}})$ denotes the signature of
$\overline{\mathcal{M}}$. Now we denote by $B(\mathcal{M})$ the
set of all restrictions of $\overline{\mathcal{M}}$ preserving the
universe $M$. Clearly, all these restrictions are regular, too.
Since each signature relation of $\overline{\mathcal{M}}$ can be
independently preserved or removed there are
$2^{\left|\Sigma(\overline{\mathcal{M}})\right|}=2^{{\rm
max}\{2^{|M|},\omega\}}$ possibilities for these restrictions,
i.e. $|B(\mathcal{M})|=2^{{\rm max}\{2^{|M|},\omega\}}$. Since
each element of $B(\mathcal{M})$ is uniquely defined by a subset
$\Sigma\subseteq\Sigma(\overline{\mathcal{M}})$, the set-theoretic
operations on the Boolean
$\mathcal{P}(\Sigma(\overline{\mathcal{M}}))$, forming its Cantor
algebra, induce the {\em regular} atomic Boolean algebra
$\mathcal{B}(\mathcal{M})$ on a universe $B(\mathcal{M})$, with
the greatest element $\overline{\mathcal{M}}$ having a complete,
i.e. maximal signature, the least element $\mathcal{M}_0$ having
the empty signature, and $|\Sigma(\overline{\mathcal{M}})|$ atoms
each of which has exactly one signature symbol. Here unions
$\mathcal{N}_1\cup\mathcal{N}_2$ and intersections
$\mathcal{N}_1\cap\mathcal{N}_2$, for
$\mathcal{N}_1,\mathcal{N}_2\in B(\mathcal{M})$, preserve the
universe $M$ and consists of unions of their signature relations,
common signature symbols, respectively. The unions can be
considered both as {\em combinations} \cite{cs} and {\em fusions}
\cite{HaHi}, in a broad sense, of structures.

The algebra $\mathcal{B}(\mathcal{M})$ includes possibilities of
regular structures on the universe $M$, and contains all
restrictions of the structure $\overline{\mathcal{M}}$. This
algebra is uniquely defined up to the names of signature symbols,
and it reflects a hierarchy of regular expansions and restrictions
of structures with respect to the names of relational symbols.
Moreover, the following theorem gives a characterization of
existence of an isomorphism between regular Boolean algebras.

\begin{theorem}\label{mod_BA}
For any regular structures $\mathcal{M}$ and $\mathcal{N}$ the
following conditions are equivalent:

$(1)$ $\mathcal{B}(\mathcal{M})\simeq\mathcal{B}(\mathcal{N})$,

$(2)$ there is a bijection between sets of atoms for
$\mathcal{B}(\mathcal{M})$ and $\mathcal{B}(\mathcal{N})$,

$(3)$ ${\rm max}\{|M|,\omega\}={\rm max}\{|N|,\omega\}$.
\end{theorem}

Proof. $(1)\Rightarrow(2)$. Let
$\mathcal{B}(\mathcal{M})\simeq\mathcal{B}(\mathcal{N})$ witnessed
by an isomorphism $f$. Then the restriction of $f$ to the set of
atoms produces a bijection between sets of atoms for
$\mathcal{B}(\mathcal{M})$ and $\mathcal{B}(\mathcal{N})$.

$(2)\Rightarrow(1)$. Since $\mathcal{B}(\mathcal{M})$ and
$\mathcal{B}(\mathcal{N})$ are uniquely defined by their atoms and
there is a bijection $f_a$ between sets of atoms for
$\mathcal{B}(\mathcal{M})$ and $\mathcal{B}(\mathcal{N})$, this
bijection $f_a$ has a unique extension till an isomorphism between
$\mathcal{B}(\mathcal{M})$ and $\mathcal{B}(\mathcal{N})$.

$(2)\Rightarrow(3)$. Since each atom for
$\mathcal{B}(\mathcal{M})$ and for $\mathcal{B}(\mathcal{N})$ is
uniquely defined by a signature symbol in
$\Sigma(\overline{\mathcal{M}})$ and in
$\Sigma(\overline{\mathcal{N}})$, respectively, the bijection
between sets of atoms for $\mathcal{B}(\mathcal{M})$ and
$\mathcal{B}(\mathcal{N})$ implies ${\rm max}\{|M|,\omega\}={\rm
max}\{|N|,\omega\}$.

$(3)\Rightarrow(2)$. Since the equality ${\rm
max}\{|M|,\omega\}={\rm max}\{|N|,\omega\}$ implies
$\left|\Sigma(\overline{\mathcal{M}})\right|=\left|\Sigma(\overline{\mathcal{N}})\right|$
in view of (\ref{eq_card_Sigma}), there is a bijection between
sets of atoms for $\mathcal{B}(\mathcal{M})$ and
$\mathcal{B}(\mathcal{N})$.
\endproof

\medskip
Theorem \ref{mod_BA} immediately implies:

\begin{corollary}\label{cor_BA}
For any at most countable regular structures $\mathcal{M}$ and
$\mathcal{N}$ their Boolean algebras $\mathcal{B}(\mathcal{M})$
and $\mathcal{B}(\mathcal{N})$ are isomorphic.
\end{corollary}

Since each algebra $\mathcal{B}(\mathcal{M})$ consists of
structures with pairwise distinct signatures all these structures
$\mathcal{M}'$ have pairwise distinct theories ${\rm
Th}(\mathcal{M}')$. Thus all structures in
$\mathcal{B}(\mathcal{M})$ can be replaced by their theories
obtaining an isomorphic Boolean algebra $\mathcal{B}({\rm
Th}(\mathcal{M}))$, which is defined by a partial order $\leq$
with $T_1\leq T_2$ iff a theory $T_2\in B({\rm Th}(\mathcal{M}))$
is an expansion of a theory $T_2\in B({\rm Th}(\mathcal{M}))$.
Here the algebra $\mathcal{B}(\mathcal{M})$ is composed by
semantic objects $\mathcal{M}'$ whereas $\mathcal{B}({\rm
Th}(\mathcal{M}))$ consists of syntactic ones: ${\rm
Th}(\mathcal{M}')$.

Notice that $\mathcal{B}({\rm Th}(\mathcal{M}))$ depends on the
cardinality $\lambda=|M|$ and, since models of their theory are
uniquely defined up to an isomorphism iff these models are finite,
one can reconstruct elements of $\mathcal{B}(\mathcal{M})$ by
elements of $\mathcal{B}({\rm Th}(\mathcal{M}))$, up to the
bijection of the universe $M$, iff $M$ is finite. So we denote
$\mathcal{B}({\rm Th}(\mathcal{M}))$ by $\mathcal{B}_\lambda(T)$,
where $T={\rm Th}(\mathcal{M})$ and $\lambda=|\mathcal{M}|$. Here
a finite cardinality $\lambda$ can be replaced by countable one by
Theorem \ref{mod_BA}.

In view of Theorem \ref{mod_BA} we have the following:

\begin{theorem}\label{th_th_BA}
For any regular theories $T_1$ and $T_2$ the following conditions
are equivalent:

$(1)$ $\mathcal{B}_\lambda(T_1)\simeq\mathcal{B}_\mu(T_2)$,

$(2)$ there is a bijection between sets of atoms for
$\mathcal{B}_\lambda(T_1)$ and $\mathcal{B}_\mu(T_2)$,

$(3)$ ${\rm max}\{\lambda,\omega\}={\rm max}\{\mu,\omega\}$.
\end{theorem}

\begin{corollary}\label{cor_th_BA}
For any countable regular theories $T_1$ and $T_2$ their Boolean
algebras $\mathcal{B}_\omega(T_1)$ and $\mathcal{B}_\omega(T_2)$
are isomorphic.
\end{corollary}

\begin{remark}\label{rem_dereg}\rm
The procedures of (de)regularizations can be naturally transformed
from structures to their theories. So by the definition any
complete theory is a deregularization of a theory $T$ in
appropriate $\mathcal{B}_\lambda(T)$.
\end{remark}

\section{Kinds of properties with respect to Boolean algebras}

Recall that for a lattice $\mathcal{L}$ and its element $a$ the
{\em upper cone} $\triangledown(a)=\triangledown_a$ consists of
all elements $b$ in $L$ with $a\leq b$, and the {\em lower cone}
$\triangle(a)=\triangle_a$ consists of all elements $b$ in $L$
with $b\leq a$.

Take a Boolean algebra $\mathcal{B}(\mathcal{M})$ and a property
$P\subseteq B(\mathcal{M})$. The property $P$ is {\em closed under
expansions} in $\mathcal{B}(\mathcal{M})$, or with respect to
$\mathcal{B}(\mathcal{M})$, or briefly the {\em
$E_{\mathcal{B}(\mathcal{M})}$-property}, if all expansions, in
$\mathcal{B}(\mathcal{M})$, of any structure $\mathcal{N}\in P$
belong to $P$. Accordingly, the property $P$ is {\em closed under
restrictions}  in $\mathcal{B}(\mathcal{M})$, or with respect to
$\mathcal{B}(\mathcal{M})$, or briefly the {\em
$R_{\mathcal{B}(\mathcal{M})}$-property}, if all restrictions, in
$\mathcal{B}(\mathcal{M})$, of any structure $\mathcal{N}\in P$
belong to $P$.

A property $P$ on the class of regular structures is {\em closed
under expansions}, or briefly the {\em $E$-property}, if all
expansions of any structure $\mathcal{N}\in P$ belong to $P$.
Accordingly, the property $P$ is {\em closed under restrictions},
or briefly the {\em $R$--property}, if all restrictions of any
structure $\mathcal{N}\in P$ belong to $P$.

\begin{remark}\label{rem_exp_restr}\rm
If a property $P\subseteq B(\mathcal{M})$ is closed both under
expansions and restrictions in $B(\mathcal{M})$ then either
$P=\emptyset$ or $P=B(\mathcal{M})$.
\end{remark}

\begin{remark}\label{rem_exp_rest}\rm Using regularizations
and deregularizations the notions of $E$-property and $R$-property
can be naturally spread for arbitrary structures and for theories,
not necessary regular.
\end{remark}

\begin{remark}\label{rem_prop_with}\rm A series of non-elementary
properties separates the closeness under expansions (restrictions)
with respect to given Boolean algebras and in general. For
instance, if $\mathcal{M}$ is a countable structure without
constant symbols then it is expansible by countably many constants
and does not expansible in $\mathcal{B}(\mathcal{M})$ by
uncountably many constants whereas $\mathcal{M}$ can be
arbitrarily many repeated constants.
\end{remark}

By the definition we have the following:

\begin{proposition}\label{prop_cone}
For any Boolean algebra $\mathcal{B}(\mathcal{M})$ and a property
$P\subseteq B(\mathcal{M})$, $P$ is closed under expansions {\rm
(}restrictions{\rm )} in $\mathcal{B}(\mathcal{M})$ iff for any
$\mathcal{N}\in P$, $\triangledown_\mathcal{N}\subseteq P$ {\rm
(}{\rm $\triangle_\mathcal{N}\subseteq P$)}, i.e.
$P=\bigcup\limits_{\mathcal{N}\in P}\triangledown_\mathcal{N}$
$\left(P=\bigcup\limits_{\mathcal{N}\in
P}\triangle_\mathcal{N}\right)$.
\end{proposition}

\begin{corollary}\label{cor_cone}
Let $P$ be a property on the class of structures closed under
regularizations and deregularizations. Then $P$ is closed under
expansions {\rm (}restrictions{\rm )} iff for any regular
$\mathcal{M}\in P$, $\triangledown_\mathcal{M}\subseteq P$ {\rm
(}{\rm $\triangle_\mathcal{M}\subseteq P$)} in any
$\mathcal{B}(\mathcal{M})$, i.e. $P$ is the deregularization of
$\bigcup\limits_{\mathcal{M}\in P}\triangledown_\mathcal{M}$
$\left(\bigcup\limits_{\mathcal{M}\in
P}\triangle_\mathcal{M}\right)$.
\end{corollary}

\begin{remark}\label{rem_tr_nab}\rm
Since elements of $\triangle_\mathcal{N}$ are obtained from
$\mathcal{N}$ just loosing some signature symbols of $\mathcal{N}$
then $\triangle_\mathcal{N}$ does not depend on choice of Boolean
algebra $\mathcal{B}(\mathcal{M})$ containing $\mathcal{N}$.
Similarly $\triangledown_\mathcal{N}$ does not depend on choice of
Boolean algebra $\mathcal{B}(\mathcal{M})$ containing
$\mathcal{N}$, since it is obtained from $\mathcal{N}$ by adding
arbitrary new signature symbols for relations preserving the
regularity. So the assertion of Corollary \ref{cor_cone} differs
from Proposition \ref{prop_cone} just by possibilities of
deregularization.
\end{remark}

By Proposition \ref{prop_cone} any complement $\overline{P}$ of a
property $P\in B(\mathcal{M})$, which is closed under expansions
(restrictions) in $\mathcal{B}(\mathcal{M})$, is represented by
some
$\bigcap\limits_{\mathcal{N}}\overline{\triangledown_\mathcal{N}}$
$\left(\bigcap\limits_{\mathcal{N}}\overline{\triangle_\mathcal{N}}\right)$.
At the same time, by the duality principle, it is closed under
restrictions (expansions). Indeed, let $\overline{P}$ be not
closed under restrictions (expansions). Then
$\mathcal{B}(\mathcal{M})$ has elements $\mathcal{N}_1$ and
$\mathcal{N}_2$ with $\mathcal{N}_1\in\overline{P}$,
$\mathcal{N}_2\in P$, $\mathcal{N}_2\in\triangle(\mathcal{N}_1)$
($\mathcal{N}_2\in\triangledown(\mathcal{N}_1)$). But then
$\mathcal{N}_1\in\triangledown(\mathcal{N}_2)$
($\mathcal{N}_1\in\triangle(\mathcal{N}_2)$) contradicting the
conditions that $P$ is closed under expansions (restrictions).

Thus, applying Proposition \ref{prop_cone} we have the following:

\begin{theorem}\label{th_property} For any property $P\subseteq
B(\mathcal{M})$ the following conditions are equivalent:

$(1)$ $P$ is closed under expansions {\rm (}restrictions{\rm )} in
$\mathcal{B}(\mathcal{M})$;

$(2)$ $P=\bigcup\limits_{\mathcal{N}\in
P}\triangledown_\mathcal{N}$
$\left(P=\bigcup\limits_{\mathcal{N}\in
P}\triangle_\mathcal{N}\right)$;

$(3)$ $\overline{P}=\bigcup\limits_{\mathcal{N}\in
\overline{P}}\triangle_\mathcal{N}$
$\left(\overline{P}=\bigcup\limits_{\mathcal{N}\in
\overline{P}}\triangledown_\mathcal{N}\right)$.
\end{theorem}

Transforming structures to their theories we obtain the following:

\begin{corollary}\label{cor_property} For any property $P\subseteq
B_\lambda(T)$ the following conditions are equivalent:

$(1)$ $P$ is closed under expansions {\rm (}restrictions{\rm )} in
$\mathcal{B}_\lambda(T)$;

$(2)$ $P=\bigcup\limits_{T'\in P}\triangledown_{T'}$
$\left(P=\bigcup\limits_{T'\in P}\triangle_{T'}\right)$;

$(3)$ $\overline{P}=\bigcup\limits_{T'\in
\overline{P}}\triangle_{T'}$
$\left(\overline{P}=\bigcup\limits_{Y'\in
\overline{P}}\triangledown_{T'}\right)$.
\end{corollary}

\section{Properties forming lattices}

{\bf Definition.} A property $P\subseteq B(\mathcal{M})$ is called
a {\em lattice property}, or {\em $L$-property}, if $P$ is the
universe of a sublattice of the lattice restriction of
$\mathcal{B}(\mathcal{M})$.

\medskip
Since $B(\mathcal{M})$ is distributive, any $L$-property is
distributive, too.

By the definition a nonempty property $P\subseteq B(\mathcal{M})$
is a $L$-property iff for any $\mathcal{N}_1,\mathcal{N}_2\in P$,
$\mathcal{N}_1\cap\mathcal{N}_2\in P$ and
$\mathcal{N}_1\cup\mathcal{N}_2\in P$. Besides lower cones
$\triangle_\mathcal{N}$ and upper cones
$\triangledown_\mathcal{N}$ are closed under intersections and
under unions, therefore they are $L$-properties. Since
$L$-properties preserve the distributivity, we have the following:

\begin{proposition}\label{prop-lat}
If $\emptyset\ne P\subseteq B(\mathcal{M})$ is closed under lower
{\rm (}upper{\rm )} cones then the following conditions are
equivalent:

$(1)$ $P$ is a $L$-property,

$(2)$ $P$ is a distributive $L$-property,

$(3)$ $P$ is closed under unions {\rm (}intersections{\rm )}.
\end{proposition}

\begin{corollary}\label{cor_lat_fin}
If $\emptyset\ne P\subseteq B(\mathcal{M})$ is closed under lower
{\rm (}upper{\rm )} cones and it is a $L$-property then the
following conditions are equivalent:

$(1)$ $P$ is a finite union of lower {\rm (}upper{\rm )} cones;

$(2)$ $P$ is a lower {\rm (}upper{\rm )} cone.
\end{corollary}

Proof. Clearly, a lower {\rm (}upper{\rm )} cone is a finite union
of lower {\rm (}upper{\rm )} cones, i.e. $(2)\Rightarrow(1)$
obviously holds. So it suffices to show, using (1) and the
$L$-property for $P$, that two lower {\rm (}upper{\rm )} cones are
reduced to one lower {\rm (}upper{\rm )} cone. Indeed, if
$\triangle(\mathcal{N}_1)\subseteq P$ and
$\triangle(\mathcal{N}_2)\subseteq P$ then $\mathcal{N}_1\cup
\mathcal{N}_2\in P$ by Proposition \ref{prop-lat}, with
$\triangle(\mathcal{N}_1)\cup\triangle(\mathcal{N}_2)\subseteq\triangle(\mathcal{N}_1\cup
\mathcal{N}_2)\subseteq P$. Similarly, if
$\triangledown(\mathcal{N}_1)\subseteq P$ and
$\triangledown(\mathcal{N}_2)\subseteq P$ then $\mathcal{N}_1\cup
\mathcal{N}_2\in P$ by Proposition \ref{prop-lat}, with
$\triangledown(\mathcal{N}_1)\cup\triangledown(\mathcal{N}_2)\subseteq\triangledown(\mathcal{N}_1\cup
\mathcal{N}_2)\subseteq P$ completing the proof of
$(1)\Rightarrow(2)$. \endproof

\begin{remark}\label{rem_L-prop}\rm
Clearly, there are many $L$-properties which are not represented
as unions of cones. For instance, closures in
$\overline{\mathcal{M}}$ with respect to unions and intersections
of sets of finitely many structures
$\mathcal{N}_1,\ldots,\mathcal{N}_k\in B(\mathcal{M})$ form finite
lattices, moreover, Boolean algebras, as any finitely generated,
i.e. finite lattices here. If signatures
$\Sigma_1,\ldots,\Sigma_k$ of
$\mathcal{N}_1,\ldots,\mathcal{N}_k$, respectively, are infinite
and co-infinite and their unions and intersections are infinite
and co-infinite, too, then these Boolean algebras do not contain
cones at all.

At the same time, if signatures $\Sigma_1,\ldots,\Sigma_k$ are
(co)finite then they can be extended by finitely many ones forming
cones of the form $\triangle_\mathcal{N}$ (respectively,
$\triangledown_\mathcal{N}$), where $\mathcal{N}$ has the
signature $\Sigma_1\cup\ldots\cup\Sigma_k$
($\Sigma_1\cap\ldots\cap\Sigma_k$).
\end{remark}

\section{Structures and properties closed under permutations}

{\bf Definition.} We say that a property $P\subseteq
B(\mathcal{M})$ is {\em closed} or {\em invariant under
permutations} if for any permutation $f$ on $M$ and for any
$\mathcal{N}\in B(\mathcal{M})$ the image $f(\mathcal{N})$, with
signature relations $f(Q)=\{f(\overline{a})\mid \overline{a}\in
Q\}$ for all given signature relations $Q$ on $\mathcal{N}$,
belongs to $P$. Here the names $Q$ for signature predicates are
preserved with respect to $f$ iff $f(Q)=Q$.

A structure $\mathcal{N}$ is {\em closed} or {\em invariant under
permutations} if the property $P=\{\mathcal{N}\}$ is closed under
permutations, i.e. each permutation of the universe of
$\mathcal{N}$ is a permutation of all signature relations.

\begin{example}\label{ex_struct_perm}\rm If $\mathcal{N}$ consists of
all signature singletons $\{a\}$, $a\in N$, then $\mathcal{N}$ is
closed under permutations. Moreover, if $\mathcal{N}$ consists of
all $n$-element unary relations, where $n$ is some natural number,
then $\mathcal{N}$ is closed under permutations, too. We observe
the same effect for relations with $n$-element complements. In
particular, if a structure $\mathcal{N}$ consists of either empty
or complete signature relations then $\mathcal{N}$ is closed under
permutations.
\end{example}

\begin{proposition}\label{prop_clos_perm}
If $P=\triangle_\mathcal{M}$ {\rm
(}$P=\triangledown_\mathcal{M}${\rm )} then the following
conditions are equivalent:

$(1)$ $P$ is closed under permutations;

$(2)$ $\mathcal{M}$ is closed under permutations.
\end{proposition}

Proof. $(1)\Rightarrow(2)$. Since $P=\triangle_\mathcal{M}$ {\rm
(}$P=\triangledown_\mathcal{M}${\rm )} and  $P$ is closed under
permutations, any permutation $f$ on $M$ preserves the relation
``to be a restriction (an expansion)''. Then more rich (poor)
structures with respect to the signature are transformed to more
rich (poor). In particular, the most rich (poor) structure
$\mathcal{M}$ in $P$ is transformed to itself, i.e. $\mathcal{M}$
is closed under permutations.

$(2)\Rightarrow(1)$. If $\mathcal{M}$ is closed under permutations
then any permutation $f$ of $M$ moves any subsignature
(supersignature) in $\mathcal{M}$ again to a subsignature
(supersignature) in $\mathcal{M}$. Then
$f(\triangle_\mathcal{M})\subseteq \triangle_\mathcal{M}$
($f(\triangledown_\mathcal{M})\subseteq
\triangledown_\mathcal{M}$). Since $f^{-1}$ is again a
permutation, we have $f(\triangle_\mathcal{M})\supseteq
\triangle_\mathcal{M}$ ($f(\triangledown_\mathcal{M})\supseteq
\triangledown_\mathcal{M}$). Thus $P$ is preserved under
permutations.
\endproof

\begin{remark}\label{rem_perm}\rm
If a property $P$ on a family $B(\mathcal{M})$ of structures is
defined by interactions of definable sets and it is preserved
under replacements of universes then $P$ is closed under
permutations. We observe the same effect for classes of all models
of theories in a given family. These classes are preserved by
isomorphisms, in particular, by permutations transforming
signature relations.
\end{remark}

\begin{example}\label{ex+perm_triangle}\rm
Let $\mathcal{M}_n$ be a structure consisting of all predicate
relations with $n$-element complements, say of given family of
arities. Then both $\mathcal{M}_n$ and $\triangle(\mathcal{M}_n)$
are invariant under permutations, with respect to any Boolean
algebra $\mathcal{B}(\mathcal{N})$, containing $\mathcal{M}_n$.

The same effect holds in expansions (restrictions) of
$\mathcal{M}_n$ by (to) empty and/or complete predicates. We
always have this invariance for one-element regular structures.
\end{example}

\section{The families of $\omega$-categorical and Ehrenfeucht structures in $\mathcal{B}(\mathcal{M})$}

Let $P=P_{\omega\mbox{\rm\footnotesize -cat}}\subseteq
B(\mathcal{M})$ be the property of $\omega$-categoricity. Recall
that it is characterized, in view of Ryll-Nardzewski theorem, by
finitely many $n$-types over the empty set, for any natural $n$.
Then this property $P_{\omega\mbox{\rm\footnotesize -cat}}$ is
invariant under permutations.

Since restrictions of theories preserve the finiteness of number
of types, $P$ is closed under lower cones: if $\mathcal{N}\in P$
then $\triangle_\mathcal{N}\subseteq P$. Thus we have
\begin{equation}\label{eq_inter}P=\bigcup\limits_{\mathcal{N}\in
B(\mathcal{M})}\triangle_\mathcal{N}.\end{equation}

The following two examples show that $P$ is not closed under
fusions.

\begin{example}\label{ex_cat1}\rm
Let $\mathcal{N}_i=\langle M;Q_i\rangle$, $i=1,2$, be two infinite
graphs whose each connected component consists of unique edge.
Clearly each $\mathcal{N}_i$ is $\omega$-categorical with finitely
many $n$-types and these types are defined by families of formulae
of forms $x\approx x$, $(x\approx y)^\delta$, $Q^\delta(x,y)$,
$\delta\in\{0,1\}$. The structure
$\mathcal{N}_1\cup\mathcal{N}_2=\langle M;Q_1,Q_2\rangle$ consists
of vertices each of which is incident to exactly two edges, of
distinct colors $Q_1$ and $Q_2$. These edges belong to either
cycles of even lengths or to infinite chains, depending on
relationship of the relations $Q_1$ and $Q_2$. In a case of
presence of unbounded lengths of cycles or an infinite chain the
fusion $\mathcal{N}_1\cup\mathcal{N}_2$ looses the
$\omega$-categoricity.
\end{example}

\begin{example}\label{ex_cat2}\rm
Let $\mathcal{N}_1=\langle M;\leq\rangle$ be a countable densely
ordered set, $\mathcal{N}_2=\langle M;R\rangle$ be a structure
with a unary predicate $R$. Obviously, both these structures are
$\omega$-categorical. Their fusion
$\mathcal{N}_1\cup\mathcal{N}_2=\langle M;\leq,R\rangle$ can
preserve or loose the $\omega$-categoricity depending on
$\leq$-ordering of elements in $R$. Indeed, if $R$ is a finite
union of $\leq$-convex sets the structure is again
$\omega$-categorical with densely ordered convex parts. And if $R$
has an infinite $\leq$-discretely ordered part then
$\mathcal{N}_1\cup\mathcal{N}_2$ is not $\omega$-categorical.

Notice that here $\mathcal{N}_1\cup\mathcal{N}_2$ can be
Ehrenfeucht, i.e. with finitely many $(>1)$ pairwise
non-isomorphic elementary equivalent countable structures, if $R$
has finitely many $(\geq 1)$ $\leq$-accumulation points in some
saturated elementary extension of
$\mathcal{N}_1\cup\mathcal{N}_2$, and continuum many pairwise
non-isomorphic elementary equivalent countable structures, if $R$
has infinitely many $\leq$-accumulation points in some saturated
elementary extension of $\mathcal{N}_1\cup\mathcal{N}_2$.
\end{example}

In view of Examples \ref{ex_cat1} and \ref{ex_cat2} the family $P$
does not form a lattice. Besides, it does not have maximal
elements since each $\omega$-categorical structure has an
$\omega$-categorical proper expansion by new signature singleton.
Thus $P$ is formed as a union of infinite increasing chains of
lower cones.

Summarizing the arguments above we have the following:

\begin{theorem}\label{th_count_cat}
The family $P_{\omega\mbox{\rm\footnotesize -cat}}\subseteq
B(\mathcal{M})$ of countably categorical structures is represented
as the union of lower cones of all its elements and all these
elements are not maximal. This family is closed under permutations
and not closed under unions.
\end{theorem}

\begin{remark}\label{rem_Ehr}\rm
In \cite{CCMCT}, a machinery described allowing to obtain
Ehrenfeucht expansions of some given non-Ehrenfeucht theories,
then loosing the Ehrenfeucht property by a suitable expansion,
returning it, etc. Thus there is a sequence of expansions of
structures in $\mathcal{B}(\mathcal{M})$, with countable
$\mathcal{M}$, alternating the Ehrenfeucht property and its
complement.

Remind that the Ehrenfeucht property can fail under inessential
expansions and restrictions, i.e. under expansions and
restrictions by constants \cite{Om}. It means that the Ehrenfeucht
property is rather different in $\mathcal{B}(\mathcal{M})$ than
the property of $\omega$-categoricity.
\end{remark}

Let $P=P_{{\rm\footnotesize Ehr}}\subseteq B(\mathcal{M})$ be the
property of Ehrenfeuchtness, where $M$ is countable. Clearly, $P$
is closed under permutations of structures, and both the least and
the greatest elements of $\mathcal{B}(\mathcal{M})$ are not
Ehrenfeucht. Therefore $P$ does not contain cones at all. At the
same time the property $P$ has atomic elements in the Boolean
algebra $\mathcal{B}(\mathcal{M})$, i.e. Ehrenfeucht structures
with exactly one signature symbol.

For instance, the classical Ehrenfeucht theory $T={\rm
Th}(\langle\mathbb Q;\leq,c_n\rangle_{i\in\omega})$,
$c_n<c_{n+1}$, $n\in\omega$, \cite{CCMCT}, with exactly three
countable models can be replaced by the theory $T'$ of the binary
structure
$$\langle(\mathbb{Q}\setminus\{c_n\mid
n\in\omega\}\times\omega)\cup\{(c_n,k)\mid k<n,n\in\omega\};$$
$$\{((q,0),(r,0))\mid q,r\in\mathbb Q, q\leq
r\}\cup\{((q,0),(q,k))\mid q\in\mathbb Q, k>0\}\rangle.$$ Here
finitely many arcs $((q,0),(q,k))$, where $q=c_n$, code the
position of the element $c_n$ among other elements in the linearly
ordered part $\{(q,0)\mid q\in\mathbb Q\}$ and produce exactly
three non-isomorphic possibilities for countable models of $T'$.

Summarizing the arguments above we have the following:

\begin{theorem}\label{th_Ehr}
Any Boolean algebra $\mathcal{B}(\mathcal{M})$ with a countable
universe $M$ contains structures with the property
$P_{{\rm\footnotesize Ehr}}$ without the least and the greatest
elements of $\mathcal{B}(\mathcal{M})$. This property is closed
under permutations and can fail under restrictions and expansions.
There are infinite chains alternating the Ehrenfeuchtness and the
complement of this property. There are atomic structures
$\mathcal{N}\in B(\mathcal{M})$ belonging to $P_{{\rm\footnotesize
Ehr}}$.
\end{theorem}

\section{The family of strongly minimal structures in $\mathcal{B}(\mathcal{M})$}

Recall \cite{BaLa} that a structure $\mathcal{N}$ is called {\em
strongly minimal} if for any $\mathcal{N}'\equiv\mathcal{N}$ and
any formula $\varphi(x,\overline{a})$ in the language of
$\mathcal{N}$ with parameters $\overline{a}\in N'$ the set
$\varphi(\mathcal{N}',\overline{a})=\{b\mid\mathcal{N}'\models\varphi(b,\overline{a})\}$
is either finite or cofinite in $N$.

A theory $T$ is called {\em strongly minimal} if $T={\rm
Th}(\mathcal{N})$ for a strongly minimal structure $\mathcal{N}$.

Let $P=P_{{\rm\footnotesize sm}}\subseteq B(\mathcal{M})$ be the
property of strong minimality.

\begin{remark}\label{rem_sm}\rm By the definition the property $P$ is closed
under permutations and intersections of structures, via
intersections of signatures, and under naming of new elements by
unary singletons. Moreover, unary predicates of finite
cardinalities can be added preserving the regularity and the
strong minimality, whereas some equivalence relations can be added
with finite bounded equivalence classes only, for instance, with
unique non-one-element class, since unbounded cardinalities imply
infinite ones which is forbidden for strongly minimal structures.
The latter condition implies that $P$ does not have maximal
elements in $\mathcal{B}(\mathcal{M})$ allowing to add new
equivalence classes with $n+1$ instead of the bound $n$.

Thus $P$ forms a lower semilattice represented by the equality
(\ref{eq_inter}). At the same time $P$ admits unions of chains of
expansions, where these expansions are obtained by unary
predicates with the finiteness condition and by equivalence
relations with the bounded finiteness condition. Hence a natural
question arises on the existence of a regular strongly minimal
structure whose all regular expansions are not strongly minimal.
Below we give an answer to this question.
\end{remark}

Now we argue to show that any fusion of strongly minimal
structures $\mathcal{N}_1,\mathcal{N}_2\in B(\mathcal{M})$ is
again strongly minimal.

\begin{theorem}\label{th_fus_sm}
If structures $\mathcal{N}_1,\mathcal{N}_2\in B(\mathcal{M})$ are
strongly minimal then $\mathcal{N}_1\cup\mathcal{N}_2$ is strongly
minimal, too.
\end{theorem}

Proof. By the conjecture any formula $\varphi(x,\overline{a})$ of
a theory $T_i={\rm Th}(\mathcal{N}_i)$, $i\in\{1,2\}$, with
parameters $\overline{a}\in M$, has either finitely many or
cofinitely many solutions. Now we take a formula
$\psi(x,\overline{b})$ of the signature
$\Sigma(\mathcal{N}_1\cup\mathcal{N}_2)$, with parameters
$\overline{b}\in M$. Assume that the set
$Z=\psi(\mathcal{N}_1\cup\mathcal{N}_2,\overline{b})$ is both
infinite and co-infinite. Since $Z$ can not be obtained by Boolean
combinations of finite and cofinite sets, we may assume that
$\psi(x,\overline{b})$ has the form $\exists
y\chi(x,y,\overline{b})$ and this formula has the minimal length
among formulae with that property. Therefore
$\chi(x,y,\overline{b})$ witnesses the strong minimality with
finite or cofinite sets of solutions for $\chi(x,c,\overline{b})$
and $\chi(c,y,\overline{b})$, $c\in M$. It implies that
$\chi(x,c,\overline{b})$ has finitely many solutions for any $c\in
M$, since otherwise, if some $d$ produces cofinite
$\chi(\mathcal{N}_1\cup\mathcal{N}_2,d,\overline{b})$, then by
$$\psi(\mathcal{N}_1\cup\mathcal{N}_2,\overline{b})=\bigcup\limits_{d'\in M}\chi(\mathcal{N}_1\cup\mathcal{N}_2,d',\overline{b})\supseteq
\chi(\mathcal{N}_1\cup\mathcal{N}_2,d,\overline{b})$$ the set
$\psi(\mathcal{N}_1\cup\mathcal{N}_2,\overline{b})$ is cofinite,
too, contradicting the assumption of its co-infinity. Similarly
for each $e\in M$ the set
$\chi(e,\mathcal{N}_1\cup\mathcal{N}_2,\overline{b})$ is finite.
Then there are cofinitely many elements $d,e\in M$ with nonempty
$\chi(\mathcal{N}_1\cup\mathcal{N}_2,d,\overline{b})$ and
$\chi(e,\mathcal{N}_1\cup\mathcal{N}_2,\overline{b})$ producing
the co-finiteness of
$\psi(\mathcal{N}_1\cup\mathcal{N}_2,\overline{b})$ which
contradicts the assumption. \endproof

\medskip
The arguments for the proof of Theorem \ref{th_fus_sm} show that
any finite family of regular strongly minimal structures has a
strongly minimal theory. Since the property of strong minimality
is reduced to formulae with special properties and reduced to the
family of finite signatures, any family of regular strongly
minimal structures $\mathcal{N}_i\in B(\mathcal{M})$, $i\in I$,
has a strongly minimal union, too. In particular, there is the
greatest strongly minimal structure $\mathcal{SM}$ in
$B(\mathcal{M})$, where subsignatures of $\Sigma(\mathcal{SM})$
form strongly minimal restrictions of $\mathcal{SM}$, including
the least element $\mathcal{N}_0$ of $\mathcal{B}(\mathcal{M})$
having the empty signature. Here, $P_{{\rm\footnotesize
sm}}=\triangle(\mathcal{SM})$. Moreover, any restriction
$\mathcal{N}$ of $\mathcal{SM}$ is strongly minimal, with strongly
minimal co-structure of the signature
$\Sigma(\mathcal{SM})\setminus\Sigma(\mathcal{N})$. Thus, the
family $B_{\rm sm}(\mathcal{M})=P_{{\rm\footnotesize sm}}$ of all
strongly minimal structures forms a Boolean algebra whose universe
is a proper subset of $B(\mathcal{M})$ which is equal to the lower
cone of $\mathcal{SM}$. Hence we have the following:

\begin{theorem}\label{th_ba_sm}
Any Boolean algebra $\mathcal{B}(\mathcal{M})$ contains a
distributive sublattice $\mathcal{B}_{\rm sm}(\mathcal{M})$ of all
strongly minimal structures $\mathcal{N}\in B(\mathcal{M})$. This
sublattice closed under permutations and forms a Boolean algebra
with the least element $\mathcal{N}_0$ and the greatest element
$\mathcal{SM}$ forming $\triangle(\mathcal{SM})$ which is equal to
$B_{\rm sm}(\mathcal{M})$.
\end{theorem}

\section{The family of $\omega_1$-categorical structures in $\mathcal{B}(\mathcal{M})$}

Let $P=P_{\omega_1\mbox{\rm\footnotesize -cat}}$ be the family of
infinite structures in $\mathcal{B}(\mathcal{M})$, with countable
languages and which are $\omega_1$-categorical, i.e. categorical
in some, i.e. any uncountable cardinality. It is known \cite{BaLa,
Yer, La} that a countable complete theory $T$ without finite
models is $\omega_1$-categorical iff it contains a $1$-cardinal
formula $\varphi(\overline{x},\overline{a})$ which is strongly
minimal, i.e. both
$|\varphi(\mathcal{N},\overline{a})|=|\mathcal{N}|$ for any
$\mathcal{N}\models T$ and each its definable subset, with
parameters, of $\varphi(\mathcal{N},\overline{a})$ is either
finite or cofinite.

Clearly, the property $P$ is closed under permutations, contains
the least element of $\mathcal{B}(\mathcal{M})$ iff $M$ is
infinite, and can fail under expansions. Indeed, if
$\mathcal{N}\in B(\mathcal{M})$ is an arbitrary
$\omega_1$-categorical structure having infinitely many
realizations of a non-algebraic complete type, then its expansion
dividing this set of realizations, by new unary predicate, into
two infinite parts is already not $\omega_1$-categorical.

The following example shows that restrictions of structures also
can violate the $\omega_1$-categoricity.

\begin{example}\label{ex_omega1-cat}\rm
Let $Q$ and $R$ be disjoint infinite unary predicates whose union
forms the universe of a structure $\mathcal{N}$ of the signature
$\langle Q,R,f\rangle$, where $f$ is s bijection between $Q$ and
$R$. The formulae $Q(x)$ and $R(x)$ are $1$-cardinal and strongly
minimal, i.e. $\mathcal{N}$ is $\omega_1$-categorical. At the same
time restrictions of $\mathcal{N}$ till one or two unary
predicates are already not $\omega_1$-categorical since the
formulae $Q(x)$ and $R(x)$ in these restrictions do not remain
$1$-cardinal under these restrictions.

Any expansion of $\mathcal{N}$ by an infinite and co-infinite
unary predicate $S\subset Q$ is not $\omega_1$-categorical.
\end{example}

\begin{remark}\label{rem_exp_rest2}\rm
Notice that the behavior of the property
$P_{\omega_1\mbox{\rm\footnotesize -cat}}$ with respect to its
non-preservation can be realized by families of unary predicates
and bijections responsible for the $1$-cardinality and obtaining
an infinite chain of expansions alternating the
$\omega_1$-categoricity and its negations. For instance, an
expansion of the structure $\mathcal{N}$ by the countable and
co-countable predicate $S$ is not $\omega_1$-categorical but any
expansion of $\langle\mathcal{N},S\rangle$ by a bijection $g$
between $S$ and $Q\setminus S$ restores the
$\omega_1$-categoricity. An appropriate new unary predicate again
violates the $\omega_1$-categoricity, a suitable new bijection
restores it, etc.

Since the $\omega_1$-categoricity excludes infinite linear orders
there are countable structures $\mathcal{N}'$ of finite signatures
with $\triangledown(\mathcal{N}')\cap
P_{\omega_1\mbox{\rm\footnotesize -cat}}=\emptyset$.
\end{remark}

Summarizing the arguments above we obtain the following:

\begin{theorem}\label{th_omega1-cat}
Any Boolean algebra $\mathcal{B}(\mathcal{M})$ with an infinite
universe $M$ contains structures with the property
$P_{\omega_1\mbox{\rm\footnotesize -cat}}$ including the least
element of $\mathcal{B}(\mathcal{M})$. This property is closed
under permutations and can fail under restrictions and expansions.
There are infinite chains alternating the $\omega_1$-categoricity
and the complement of this property. There are structures
$\mathcal{N}\in B(\mathcal{M})$ of finite signatures with
$\triangledown_\mathcal{N}\cap P_{\omega_1\mbox{\rm\footnotesize
-cat}}=\emptyset$.
\end{theorem}

\section{The family of stable structures in $\mathcal{B}(\mathcal{M})$}

Recall \cite{Sh} that a formula
$\varphi(\overline{x},\overline{y})$ of a theory $T$ is called
{\em stable} if there are no tuples
$\overline{a}_i,\overline{b_i}\in N$, where $i\in\omega$,
$\mathcal{N}\models T$, such that
$$\mathcal{N}\models\varphi(\overline{a}_i,\overline{b}_j)\Leftrightarrow i\leq j.$$
The theory $T$ is called {\em stable} if all its formulae are
stable. Models of a stable theory are called {\em stable}, too. Is
a formula/theory/structure is not stable, it is called {\em
unstable} or having the {\em order property}.

Following \cite{Sh} it is said that a formula
$\varphi(\overline{x},\overline{y})$ has the {\em strict order
property} if there are parameters $\overline{a}_i\in N$,
$i\in\omega$, such that the sets
$\varphi(\overline{a}_i,\mathcal{N})$, $i\in\omega$, form a
strictly descending chain with
$\varphi(\overline{a}_i,\mathcal{N})\supsetneqq
\varphi(\overline{a}_{i+1},\mathcal{N})$, $i\in\omega$.

Again following \cite{Sh} it is said that an unstable formula
$\varphi(\overline{x},\overline{y})$ has the {\em independence
property} if in every/some model $\mathcal{N}$ of $T$ there is,
for each $n\in\omega$, a family of tuples $\overline{a}_i$, $i\in
n$, such that for each of the $2^n$ subsets $X$ of $n$ there is a
tuple $\overline{b}\in N$ for which
$$\mathcal{N}\models\varphi(\overline{a}_i,\overline{b})\Leftrightarrow i\in X.$$

It is well known that a formula
$\varphi(\overline{x},\overline{y})$ is unstable iff it has the
strict order property or the independence property.

Let $P=P_{{\rm\footnotesize st}}\subseteq B(\mathcal{M})$ be the
property of stability of a structure.

\begin{remark}\label{rem_st}\rm By the definition, similarly to the strongly minimal case, the property $P$ is closed under
intersections of structures, via intersections of signatures, and
under marking of subsets of $M$ by unary predicates.

Thus $P$ forms a lower semilattice represented by the equality
(\ref{eq_inter}). At the same time $P$ admits unions of chains of
regular expansions by new unary predicates preserving the
regularity of structures.
\end{remark}

Remind \cite{HaHa} that Boolean combinations of stable formulae
stay stable. So  the violations of stability for fusions of stable
structures can be obtained via quantifiers only. The following
example illustrates this violation showing that indeed fusions of
stable structures can fail the property of stability.

\begin{example}\label{ex_st_nst}\rm
Let $R_1=\{\langle a_i,c_{ik}\mid i,k\in\omega\}$, $R_2=\{\langle
d_{lj},b_j\mid l,j\in\omega\}$ be binary relations with pairwise
distinct elements $a_i,b_i,c_{ik},d_{lj}$. We identify some
elements $c_{ik}$ and $d_{lj}$ as follows: $c_{ik}=d_{lj}
\Leftrightarrow i\leq k=l\leq j$. The obtained structure
$\mathcal{N}=\langle M;R_1,R_2\rangle$, with the universe
consisting of the elements $a_i,b_i,c_{ik},d_{lj}$ having the
identifications above only, is unstable since the formula
$\varphi(x,y)=\exists z(R_1(x,z)\wedge R_2(z,y))$ witnesses the
strict order property fails by the sequences $a_i,b_j$,
$i,j\in\omega$. At the same time the restrictions
$\mathcal{N}_i=\mathcal{N}|_{R_i}$, $i=1,2$, are stable consisting
of countably many pairwise isomorphic countable acyclic connected
components of diameter $2$ and countably many isolated elements.
Here $\mathcal{N}=\mathcal{N}_1\cup\mathcal{N}_2$.

Besides this example again confirms that $\omega$-categorical
structures $\mathcal{N}_1,\mathcal{N}_2$ can have a fusion
$\mathcal{N}_1\cup\mathcal{N}_2$ which is not
$\omega$-categorical.
\end{example}

\begin{example}\label{ex_st_nst2}\rm
Example \ref{ex_st_nst} can be easily modified to the random
bipartite graph with the relation $R$ defined by the same formula
$\varphi(x,y)$ as follows. We enumerate the family $Z$ of all
pairs of finite disjoint sets $X,Y\subset\omega$: $\nu\mbox{\rm :
}\omega\to Z$. Now for each $i\in\omega$ and $\nu(i)=\langle X,
Y\rangle$, we identify the elements $c_{ik}$ and $d_{lj}$ iff
$k=l=j$ and $j\in X$. Thus after this identification we obtain the
structure $\mathcal{N}=\langle M;R_1,R_2\rangle$ with
$R(a_i,\mathcal{N})=X$ producing the random bipartite graph having
the independence property.

The restrictions $\mathcal{N}_i=\mathcal{N}|_{R_i}$, $i=1,2$, are
again stable consisting of countably many pairwise isomorphic
countable acyclic connected components of diameter $2$ and
countably many isolated elements. Here again
$\mathcal{N}=\mathcal{N}_1\cup\mathcal{N}_2$.
\end{example}

\begin{theorem}\label{th_stab}
The family $P_{{\rm\footnotesize st}}\subseteq B(\mathcal{M})$ of
stable structures is represented as the union of lower cones of
all its elements. This family is closed under permutations and not
closed under unions and these unions can produce both the strict
order property and the independence property.
\end{theorem}

\section{Conclusion}
We considered hierarchies of properties on a Boolean algebra of
regular expansions and restrictions based on the universe of an
arbitrary regular structure. Some general properties of these
hierarchies are studied with respect to elementary theories of
structures, lower and upper cones, lattices, permutations. A
general approach is applied for several natural classes of
structures and their theories including the $\omega$-categoricity,
Ehrenfeuchtness, strong minimality, $\omega_1$-categoricity,
stability. Properties of these classes in a Boolean algebra are
described. It would be interesting to spread this approach
describing properties of further natural classes in these Boolean
algebras.

\end{document}